\documentclass{amsart}

\usepackage{fancyhdr}
\usepackage{lastpage}
\usepackage{stmaryrd,yhmath}

\newcommand{\bdis}{\begin{displaymath}}
\newcommand{\edis}{\end{displaymath}}
\newcommand{\be}{\begin{equation}}
\newcommand{\ee}{\end{equation}}
\newcommand{\mbb}{\mathbb}
\newcommand{\mcal}{\mathcal}

\newcommand{\vp}{\varphi}

\newcommand{\vth}{\vartheta}

\newcommand{\zf}{\zeta\left(\frac{1}{2}+it\right)}

\theoremstyle{definition}

\theoremstyle{remark}
\newtheorem{remark}[]{Remark}

\newtheorem*{mydef1}{{\bf Theorem}}

\newtheorem*{mydef41}{{\bf Corollary 1}}

\newtheorem*{mydef42}{{\bf Corollary 2}}

\newtheorem*{mydef43}{{\bf Corollary 3}}

\newtheorem*{mydef61}{{\bf Example 1}}

\newtheorem*{mydef62}{{\bf Example 2}}

\newtheorem*{mydef63}{{\bf Example 3}}

\numberwithin{equation}{section}

\begin{document}

\title{Jacob's ladders and properties of complete additivity and complete multiplicativity in the set of reverse
iterated integrals (energies)}

\author{Jan Moser}

\address{Department of Mathematical Analysis and Numerical Mathematics, Comenius University, Mlynska Dolina M105, 842 48 Bratislava, SLOVAKIA}

\email{jan.mozer@fmph.uniba.sk}


\begin{abstract}
New class of integral identities concerning constraints on behavior of the Riemann's zeta function on the critical line is introduced in this paper.
Namely, we have obtained new kind of $\sigma$-additivity and $\sigma$-multiplicativity in the class of reverse iterated integrals (energies).
\end{abstract}
\maketitle

\section{Introduction}

\subsection{}

Let us remind we have proved the following theorem (see \cite{3}, (2.1) -- (2.7)): for every $L_2$-orthogonal system
\bdis
\{f_n(t)\}_{n=1}^\infty,\ t\in [0,2l],\ l=o\left(\frac{T}{\ln T}\right),\ T\to\infty
\edis
there is a continuum set of $L_2$-orthogonal systems
\bdis
\begin{split}
& \{ F_n(t;T,k,l)\}_{n=1}^\infty= \\
& = \left\{ f_n(\vp_1^k(t)-T)\prod_{r=0}^{k-1}|\tilde{Z}[\vp_1^r(t)]|\right\},\
t\in [\overset{k}{T},\overset{k}{\wideparen{T+2l}}],
\end{split}
\edis
where
\bdis
\begin{split}
& \vp_1\left\{ [\overset{k}{T},\overset{k}{\wideparen{T+2l}}]\right\}=
[\overset{k-1}{T},\overset{k-1}{\wideparen{T+2l}}],\ k=1,\dots,k_0, \\
& [\overset{0}{T},\overset{0}{\wideparen{T+2l}}]=[T,T+2l],\ T\to\infty,
\end{split}
\edis
i. e. the following formula is valid
\be \label{1.1}
\begin{split}
& \int_{\overset{k}{T}}^{\overset{k}{\wideparen{T+2l}}}f_m(\vp_1^k(t)-T)f_n(\vp_1^k(t)-T)\prod_{r=0}^{k-1}\tilde{Z}^2[\vp_1^r(t)]{\rm d}t=\\
& =
\left\{\begin{array}{rcl} 0 & , & m\not=n, \\ A_n & , & m=n, \end{array} \right. , \\
& A_n=\int_0^{2l} f_n^2(t){\rm d}t.
\end{split}
\ee
Next, we have (see \cite{2}, (9.1), (9.2)) that
\bdis
\begin{split}
& \tilde{Z}^2(t)=\frac{{\rm d}\vp_1(t)}{{\rm d}t}=\frac{Z^2(t)}{2\Phi'_\vp[\vp(t)]}=
\frac{\left|\zf\right|^2}{\omega(t)},\\
& \vp_1(t)=\frac 12\vp(t) , \\
& \omega(t)=\left\{ 1+\mcal{O}\left(\frac{\ln\ln t}{\ln t}\right)\right\}\ln t,
\end{split}
\edis
and (see \cite{6}, pp. 77, 329)
\be \label{1.2}
\begin{split}
& Z(t)=e^{i\vth(t)}\zf, \\
& \vth(t)=-\frac t2\ln\pi+\text{Im}\ln\Gamma\left(\frac 14+i\frac t2\right).
\end{split}
\ee
Hence, for the classical Fourier's orthogonal system, for example,
\be \label{1.3}
\begin{split}
& \left\{ 1, \cos\frac{\pi t}{l},\sin\frac{\pi t}{l},\dots,\cos\frac{n\pi t}{l},\sin\frac{n\pi t}{l},\dots\right\}, \\
& t\in [0,2l],
\end{split}
\ee
we have the following continuum set of orthogonal systems according to (\ref{1.1})
\bdis
\begin{split}
 & \left\{ \prod_{r=0}^{k-1}\frac{\left|\zeta\left(\frac 12+i\vp_1^r(t)\right)\right|}
 {\sqrt{\omega[\vp_1^r(t)]}},\dots,\right.  \\
 & \left. \left(\prod_{r=0}^{k-1}\frac{\left|\zeta\left(\frac 12+i\vp_1^r(t)\right)\right|}
 {\sqrt{\omega[\vp_1^r(t)]}}\right)\cos\left(\frac \pi l n(\vp_1^r(t)-T)\right),\right.  \\
 & \left. \left(\prod_{r=0}^{k-1}\frac{\left|\zeta\left(\frac 12+i\vp_1^r(t)\right)\right|}
 {\sqrt{\omega[\vp_1^r(t)]}}\right)\sin\left(\frac \pi l n(\vp_1^r(t)-T)\right),\dots\right\} , \\
 & t\in [\overset{k}{T},\overset{k}{\wideparen{T+2l}}],\
 k=1,\dots,k_0, \\
 & T \geq T[\vp_1].
\end{split}
\edis

\subsection{}

We have already noticed in our paper \cite{3} that the formula (\ref{1.1}) can serve as a resource for new
integral identities in the theory of the Riemann zeta-function. For example, in the case of the first function
of the Fourier's system (\ref{1.3})
\bdis
f_1(t)=1 ,
\edis
we have the following formula (see (\ref{1.1}))
\be \label{1.4}
\int_{\overset{k}{T}}^{\overset{k}{\wideparen{T+2l}}}
\prod_{r=0}^{k-1}\tilde{Z}^2[\vp_1^r(t)]{\rm d}t=2l,\quad k=1,\dots,k_0,\quad T\to\infty.
\ee

In this paper we shall interpret the formula (\ref{1.4}) as new kind of unit operator as well as new kind
of parametric integral. Consequently, we obtain new class of integral identities for the function
\bdis
\tilde{Z}^2(t)=\frac{\left|\zf\right|^2}{\omega(t)},
\edis
i. e. (see (\ref{1.2})) the constraints on a behavior of weakly-modulated function
\bdis
\left|\zf\right|^2.
\edis
Namely, we obtain new kind of $\sigma$-additivity and also $\sigma$-multiplicativity in the class of
reversely iterated integrals (energies).

\section{Theorem}

\subsection{}

Based on our formula (\ref{1.4}) we give the following

\begin{mydef1}
 Let us denote by $G(S)$ the class of functions
 \bdis
 \begin{split}
  & g=g(u_1,\dots,u_n) , \\
  & (u_1,\dots,u_n)\in S\subset \mbb{R}^n
 \end{split}
\edis
such that
\be \label{2.1}
g\geq 0,\quad g=o\left(\frac{T}{\ln T}\right),\ T\to\infty.
\ee
Then we have the following formula
\be \label{2.2}
\begin{split}
 & \forall g\in G(S):\\
 & \int_{\overset{k}{T}}^{\overset{k}{\wideparen{T+g}}}
 \prod_{r=0}^{k-1}\tilde{Z}^2[\vp_1^r(t)]{\rm d}t=g,\quad k=1,\dots,k_0
\end{split}
\ee
for every fixed $k_0\in\mbb{N}$ and for every sufficiently big $T>0$.
\end{mydef1}

Let us remind the following properties connected with complicated structure of the formula (\ref{2.2}). First of all
we have (comp. \cite{3}, (2.5) -- (2.7)) the following
\bdis
g=o\left(\frac{T}{\ln T}\right) \ \Rightarrow
\edis
\be \label{2.3}
|[\overset{k}{T},\overset{k}{\wideparen{T+g}}]|=\overset{k}{\wideparen{T+g}}-\overset{k}{T}=o\left(\frac{T}{\ln T}\right),
\ee
\be \label{2.4}
|[\overset{k-1}{\wideparen{T+g}},\overset{k}{T}]|\sim (1-c)\pi(T);\quad \pi(T)\sim \frac{T}{\ln T},
\ee
\be \label{2.5}
[T,T+g]\prec [\overset{1}{T},\overset{1}{\wideparen{T+g}}]\prec \dots \prec [\overset{k}{T},\overset{k}{\wideparen{T+g}}]\prec \dots
\ee
where $\pi(T)$ stands for the prime-counting function and $c$ is the Euler's constant.

\begin{remark}
 Asymptotic behavior of the following disconnected set (see (\ref{2.3}) -- (\ref{2.5}), comp. \cite{3}, (2.9))
\be \label{2.6}
\Delta(T,k,g)=\bigcup_{r=0}^k [\overset{r}{T},\overset{r}{\wideparen{T+g}}]
\ee
is as follows: if $T\to\infty$ then the components of the set (\ref{2.6}) recede unboundedly each from other and all together
are receding to infinity. Hence, if $T\to\infty$, the set (\ref{2.6}) behaves as one dimensional Friedmann-Hubble expanding
universe.
\end{remark}

\begin{remark}
 For the arguments of functions in (\ref{2.2}) we have: if
 \bdis
 t\in [\overset{k}{T},\overset{k}{\wideparen{T+g}}],\quad k=1,\dots,k_0
 \edis
 then
 \bdis
 \vp_1^r(t)\in [\overset{k-r}{T},\overset{k-r}{\wideparen{T+g}}],\quad r=0,1,\dots,k,
 \edis
 (comp. \cite{3}, (2.10)).
\end{remark}

\subsection{}

Let us notice the following about the interpretation of the formula (\ref{2.2}).

\begin{remark}
New type of operator
\be \label{2.7}
\hat{H}O=\int_{\overset{k}{T}}^{\overset{k}{\wideparen{T+O}}}\prod_{r=0}^{k-1}\tilde{Z}^2[\vp_1^r(t)]{\rm d}t,\ k=1,\dots,k_0
\ee
is defined by the our formula (\ref{2.2}). This operator acts via the upper limit of the reversely iterated integral as
follows
\bdis
\begin{split}
 & \hat{H}g=
 \left(\int_{\overset{k}{T}}^{\overset{k}{\wideparen{T+O}}}\prod_{r=0}^{k-1}\tilde{Z}^2[\vp_1^r(t)]{\rm d}t\right)g= \\
 & = \int_{\overset{k}{T}}^{\overset{k}{\wideparen{T+g}}}\prod_{r=0}^{k-1}\tilde{Z}^2[\vp_1^r(t)]{\rm d}t.
\end{split}
\edis
We have, of course, that (see (\ref{2.2}))
\bdis
\hat{H}g=g,\quad \forall g\in G(S),
\edis
i.e. $\hat{H}$ is the unit operator.
\end{remark}

\begin{remark}
A kind of \emph{closed analytic cycle} on $G(S)$ is defined by our formula (\ref{2.2}). Namely, we have
\bdis
\begin{split}
 & g\longrightarrow [T,T+g]\longrightarrow [\overset{k}{T},\overset{k}{\wideparen{T+g}}] \longrightarrow \\
 & \longrightarrow \int_{\overset{k}{T}}^{\overset{k}{\wideparen{T+g}}}\prod_{r=0}^{k-1}\tilde{Z}^2[\vp_1^r(t)]{\rm d}t=g,\quad
 k=1,\dots,k_0.
\end{split}
\edis
\end{remark}

\section{Complete additivity of the reversely iterated integrals (energies). Comparison with the $\sigma$-additivity}

\subsection{}

First of all, the following corollary holds true.

\begin{mydef41}
Let
\bdis
0\leq g_l(u_1,\dots,u_n)=g(l),\ (u_1,\dots,u_n)\in S
\edis
and
\be \label{3.1}
g=\sum_{l=1}^\infty g_l\in G(S).
\ee
Then
\be \label{3.2}
\begin{split}
 & \int_{\overset{k}{T}}^{\overset{k}{\wideparen{T+\sum g(l)}}} \prod_{r=0}^{k-1} \tilde{Z}^2[\vp_1^r(t)]{\rm d}t=\\
 & = \sum_{l=1}^\infty
 \int_{\overset{k(l)}{T}}^{\overset{k(l)}{\wideparen{T+g(l)}}}
 \prod_{r=0}^{k(l)-1} \tilde{Z}^2[\vp_1^r(t)]{\rm d}t,
\end{split}
\ee
\be \label{3.3}
k,k(l)=1,\dots,k_0,\quad T\to\infty.
\ee
\end{mydef41}

\begin{remark}
We shall call the property (\ref{3.2}) as \emph{the complete additivity of the reversely iterated
integrals (energies)}.
\end{remark}

\subsection{}

Next, we give the following

\begin{mydef61}
If
\be \label{3.4}
k=17,\ l=1,2
\ee
then we obtain from (\ref{3.2}) by using the mean-value theorem (comp. \cite{3}, (4.3) -- (4.5)) that
\bdis
\begin{split}
 & \int_{\overset{17}{T}}^{\overset{17}{\wideparen{T+g(1)+g(2)}}} \prod_{r=0}^{16}
 \left|\zeta\left(\frac 12+i\vp_1^r(t)\right)\right|^2{\rm d}t\sim \\
 & \sim \ln^{16}T\cdot \int_{\overset{1}{T}}^{\overset{1}{\wideparen{T+g(1)}}}\left|\zf\right|^2{\rm d}t+ \\
 & + \ln^{10}T\cdot \int_{\overset{7}{T}}^{\overset{7}{\wideparen{T+g(2)}}}\prod_{r=0}^6
 \left|\zeta\left(\frac 12+i\vp_1^r(t)\right)\right|^2{\rm d}t,\quad T\to\infty
\end{split}
\edis
together with other formulae of the complete finite set defined by the condition (\ref{3.4}).
\end{mydef61}

\subsection{}

Furthermore, we define following planar figures (comp. \cite{4}, (3.1) -- (3.3)):
\be \label{3.5}
\begin{split}
 & P_k(T,g)=\left\{ (t,y):\ t\in [\overset{k}{T},\overset{k}{\wideparen{T+\sum g(l)}}], \
 y\in [0,\prod_{r=0}^{k-1} \tilde{Z}^2[\vp_1^r(t)]]\right\} , \\
 & k=1,\dots,k_0,
\end{split}
\ee
\be \label{3.6}
\begin{split}
 & P_{k(l)}^l(T,g(l))=\left\{ (t,y):\ t\in [\overset{k(l)}{T},\overset{k(l)}{\wideparen{T+g(l)}}], \
 y\in [0,\prod_{r=0}^{k(l)-1} \tilde{Z}^2[\vp_1^r(t)]]\right\} , \\
 & k=1,\dots,k_0, \ l=1,2,\dots
\end{split}
\ee

\begin{remark}
Consequently, we have (see (\ref{3.2})) for measures of these planar figures the following
\be \label{3.7}
m\{ P_k\}=\sum_{l=1}^\infty m\{ P_{k(l)}^l\},
\ee
where
\be \label{3.8}
\{ l_1\not=l_2 \ \& \ k(l_1)=k(l_2)\}\ \Rightarrow \ P_{k(l_1)}^{l_1}\bigcap P_{k(l_2)}^{l_2}=\emptyset .
\ee
\end{remark}

Of course, the sequence
\bdis
\{ k(l)\}_{l=1}^\infty
\edis
contains an ininite set for which
\bdis
k(l_1)=k(l_2)
\edis
(see (\ref{3.3}); $k_0$ being fixed).

\subsection{}

Let us remind that the Lebesgue measure is countable additive, or $\sigma$-additive, i.e. if
\bdis
A_1,A_2,\dots,A_n,\dots
\edis
are measurable sets and if they are pairwise disjoint
\be \label{3.9}
p\not=q \ \Rightarrow \ A_p\bigcap A_q=\emptyset ,
\ee
then
\bdis
m\{ A\}=m\left\{\bigcup_{k=1}^\infty A_k\right\}=\sum_{k=1}^\infty m\{ A_k\}.
\edis

\begin{remark}
Hence, we see (comp. (\ref{3.8}), (\ref{3.9})) that the complete additivity in our using differs from
the $\sigma$-additivity. In our case (\ref{3.7}) the Jordan's measure would be sufficient.
\end{remark}

\section{Complete multiplicativity of the reversely iterated integrals (energies)}

\subsection{}

Next, the following corollary holds true.

\begin{mydef42}
Let\footnote{Of course
\bdis
g(l_1)g(l_2) \in G(S)
\edis
does not imply
\bdis
g(l_1)\in G(S) \ \& \ g(l_2)\in G(S).
\edis}
\bdis
g_l>0;\quad g_l, \prod_{l=1}^\infty g_l\in G(S),
\edis
then
\be \label{4.1}
\begin{split}
 & \int_{\overset{k}{T}}^{\overset{k}{\wideparen{T+\prod g(l)}}}\prod_{r=0}^{k-1}\tilde{Z}^2[\vp_1^r(t)]{\rm d}t=
 \prod_{l=1}^\infty \int_{\overset{k(l)}{T}}^{\overset{k(l)}{\wideparen{T+g(l)}}}
 \prod_{r=0}^{k(l)-1}\tilde{Z}^2[\vp_1^r(t)]{\rm d}t, \\
 & k,k(l)=1,\dots,k_0,\ T\to\infty.
\end{split}
\ee
\end{mydef42}

\begin{remark}
We shall call the property (\ref{4.1}) as \emph{the complete multiplicativity of the reversely iterated
integrals (energies)}.
\end{remark}

\begin{remark}
We obtain from (\ref{3.5}), (\ref{3.6}) and (\ref{4.1}) that
\be \label{4.2}
m\{ Q_k\}=\prod_{l=1}^\infty m\{ P_{k(l)}^l\},
\ee
where in (\ref{3.5})
\bdis
P_k\left( \sum g(l) \longrightarrow \prod g(l)\right)=Q_k,
\edis
and the property (\ref{3.8}) holds true.
\end{remark}

\begin{remark}
Formulae (\ref{3.2}), (\ref{4.1}) (i.e. (\ref{3.7}), (\ref{4.2})) are not accessible by the current methods in the
theory of the Riemann zeta-function.
\end{remark}

\subsection{}

Now we give two examples.

\begin{mydef62}
If
\be \label{4.3}
k=17,\ l=1,2
\ee
then we obtain from (\ref{4.1}) by the usual way (comp. Example 1) that
\bdis
\begin{split}
 & \int_{\overset{1}{T}}^{\overset{1}{\wideparen{T+g(1)g(2)}}}\left|\zf\right|^2{\rm d}t\sim \\
 & \sim \frac{1}{\ln^{23}T}
 \int_{\overset{17}{T}}^{\overset{17}{\wideparen{T+g(1)}}}
 \prod_{r=0}^{16}\left|\zeta\left(\frac 12+i\vp_1^r(t)\right)\right|^2{\rm d}t\cdot
 \int_{\overset{7}{T}}^{\overset{7}{\wideparen{T+g(2)}}}
 \prod_{r=0}^{6}\left|\zeta\left(\frac 12+i\vp_1^r(t)\right)\right|^2{\rm d}t, \\
 & T\to\infty.
\end{split}
\edis
This formula represents the complete finite (for fixed $T,g(1),g(2)$)
set of formulae defined by the condition (\ref{4.3}).
\end{mydef62}

\begin{mydef63}
Canonical arithmetic formula
\bdis
n=p_1^{\alpha_1}p_2^{\alpha_2}\cdots p_s^{\alpha_s}\in G(S),\ n\in\mbb{N};\ \alpha_s=\alpha(s)
\edis
corresponds to the following formula
\bdis
\begin{split}
 & (\ln T)^{\sum \alpha(s)k(l)-k}
 \int_{\overset{k}{T}}^{\overset{k}{\wideparen{T+n}}}\prod_{r=0}^{k-1}
 \left|\zeta\left(\frac 12+i\vp_1^r(t)\right)\right|^2{\rm d}t\sim \\
 & \sim \prod_{l=1}^s
 \left\{ \int_{\overset{k(l)}{T}}^{\overset{k(l)}{\wideparen{T+p(l)}}}
 \prod_{r=0}^{k(l)-1}\left|\zeta\left(\frac 12+i\vp_1^r(t)\right)\right|^2{\rm d}t \right\}^{\alpha(l)},\quad
 T\to\infty.
\end{split}
\edis
\end{mydef63}

\section{On a set of asymptotically equivalent integrals connected with the length of the
Riemann's curve}

\subsection{}

Let us denote the roots of the equations
\bdis
Z(t)=0,\quad Z'(t)=0
\edis
by the symbols
\bdis
\{ \gamma\},\ \{ t_0\};\quad \gamma\not=t_0
\edis
correspondingly.

\begin{remark}
On the Riemann hypothesis the points of sequences $\{\gamma\}$ and $\{ t_0\}$ are separated each from other
(see \cite{1}, Cor. 3), i.e. in this case we have
\bdis
\gamma'<t_0<\gamma'' ,
\edis
where $\gamma',\gamma''$ are neighbouring points of the sequence $\{\gamma\}$. Of course, the value $Z(t_0)$ is
then locally extremal value of the function $Z(t)$ located at the point $t=t_0$.
\end{remark}

Next, we have proved (see \cite{5}, (1.5)) the following asymptotic formula. On the Riemann hypothesis we have
\be \label{5.1}
\begin{split}
 & \int_T^{T+H}\sqrt{1+\{ Z'(t)\}^2}{\rm d}t=2\cdot \sum_{T\leq t_0\leq T+H}|Z(t_0)|+ \\
 & + \Theta H+\mcal{O}\left( T^{\frac{A}{\ln\ln T}}\right), \\
 & \Theta=\Theta(T,H)\in (0,1),\quad H=T^{\epsilon},\ T\to\infty,
\end{split}
\ee
for every small and fixed $\epsilon>0$.

\begin{remark}
Geometric meaning of the formula (\ref{5.1}) is as follows: the length of the Riemann's curve
\bdis
y=Z(t),\quad t\in [T,T+H]
\edis
is asymptotically equal to the double of the sum of local maxima of the function
\bdis
|Z(t)|,\ t\in [T,T+H].
\edis
\end{remark}

\subsection{}

Since we have by Remark 11 that
\bdis
\sum_{T\leq t_0\leq T+H}1=\mcal{O}(H\ln T) ,
\edis
and (see \cite{6}, p. 237)
\bdis
Z(t)=\mcal{O}\left( T^{\frac{A}{\ln\ln T}}\right),
\edis
then (comp. (\ref{2.1}))
\bdis
\begin{split}
 & \sum_{T\leq t_0\leq T+H} |Z(t_0)|=
 \mcal{O}\left( H T^{\frac{A}{\ln\ln T}}\ln T\right)=
 \mcal{O}\left( T^{\epsilon+\frac{A}{\ln\ln T}}\ln T\right)= \\
 & = \mcal{O}(T^{2\epsilon})=o\left(\frac{T}{\ln T}\right).
\end{split}
\edis
Consequently, we have the following (see also the formula \cite{5}, (A.1))

\begin{mydef43}
On the Riemann hypothesis we have
\be \label{5.2}
\begin{split}
 & \int_T^{T+H}\sqrt{1+\{ Z'(t)\}^2}{\rm d}t=
 \int_{\overset{1}{T}}^{\overset{1}{\wideparen{T+H}}}
 \sqrt{1+\{ Z'_{\vp_1}[\vp_1(t)]\}^2}{\rm d}t\tilde{Z}^2(t){\rm d}t\sim \\
 & \sim \int_{\overset{k}{T}}^{\overset{k}{\wideparen{T+\sum 2|Z(t_0)|}}}\prod_{r=0}^{k-1}
 \tilde{Z}^2[\vp_1^r(t)]{\rm d}t\sim \\
 & \sim \sum_{T\leq t_0\leq T+H}
 \int_{\overset{k(t_0)}{T}}^{\overset{k(t_0)}{\wideparen{T+ 2|Z(t_0)|}}}
 \prod_{r=0}^{k(t_0)-1}\tilde{Z}^2[\vp_1^r(t)]{\rm d}t,\\
 & k,k(t_0)=1,\dots,k_0,\
 T\to\infty.
\end{split}
\ee
\end{mydef43}

\begin{remark}
The formula (\ref{5.2}) contains  asymptotic expressions of irrational integrals in
$Z'$ by rational integrals in $\tilde{Z}$.
\end{remark}

\begin{remark}
We emphasise the following simple formula
\be \label{5.3}
\begin{split}
 & \int_T^{T+H}\sqrt{1+\{ Z'(t)\}^2}{\rm d}t\sim \\
 & \sim \sum_{T\leq t_0\leq T+H}
 \int_{\overset{1}{T}}^{\overset{1}{\wideparen{T+ 2|Z(t_0)|}}}
 \tilde{Z}^2(t){\rm d}t,\quad T\to\infty,
\end{split}
\ee
as well as the formula
\bdis
\begin{split}
 & \int_T^{T+H}\sqrt{1+\{ Z'(t)\}^2}{\rm d}t\sim \\
 & \sim \frac{1}{\ln T}\sum_{T\leq t_0\leq T+H}
 \int_{\overset{1}{T}}^{\overset{1}{\wideparen{T+ 2|Z(t_0)|}}}
 \left|\zf\right|^2{\rm d}t,\quad T\to\infty.
\end{split}
\edis
\end{remark}

\thanks{I would like to thank Michal Demetrian for his help with electronic version of this paper.}

\end{document}